\documentclass[12pt,thmsa]{article}
\usepackage{sw20bams}

\input tcilatex
\begin{document}

\author{Steven R. Finch}
\title{Three Random Tangents to a Circle}
\date{January 20, 2011}
\maketitle

\begin{abstract}
Among several things, we find the side density for random triangles
circumscribing the unit circle and calculate that its median is $5.5482....$
An analogous exact computation for perimeter density remains open.
\end{abstract}

\footnotetext{
Copyright \copyright\ 2011 by Steven R. Finch. All rights reserved.}Let us
initially discuss two random tangents to the unit circle. Without loss of
generality, let one of the tangents be the vertical line passing through the
point $(-1,0)$. Let the other tangent pass through the point $(\cos (\theta
),\sin (\theta ))$, where $\theta $ is uniformly distributed on the interval 
$[0,\pi ]$. Hence it has slope $-\cot (\theta )$ and is the unique such line
touching the upper semicircle. The two lines cross at $(x,y)$, where 
\[
\left\{ 
\begin{array}{l}
x=-1, \\ 
y-\sin (\theta )=-\cot (\theta )\left( x-\cos (\theta )\right)
\end{array}
\right. 
\]
thus $h=\cot (\theta /2)$ is the (positive) height of the intersection
point. We wish to determine the probability density of $h$. Since 
\[
\frac d{d\theta }\cot \left( \frac \theta 2\right) =-\frac 12\csc \left(
\frac \theta 2\right) ^2 
\]
the density is \cite{Pa} 
\[
\left. \frac{\frac 1\pi }{\frac 12\csc \left( \frac \theta 2\right) ^2}%
\right| _{\theta =2\limfunc{arccot}(h)}=\left. \frac{\frac 2\pi }{\cot
\left( \frac \theta 2\right) ^2+1}\right| _{\theta =2\limfunc{arccot}%
(h)}=\frac 2\pi \frac 1{h^2+1} 
\]
for $h>0$, that is, the one-sided Cauchy distribution. The mean of $h$ is
infinite, as is well-known; its median is $1$.

Let us now discuss three random tangents to the unit circle, incorrectly
modeled. Take first and second lines exactly as in the preceding, and define
similarly a third line to touch the lower semicircle, independent of the
second. We study $h+k$, where $h$ is as before and $k$ is the (positive)
depth of the intersection between first and third lines. Let $\ell =h+k$.
The density of $\ell $ is the convolution \cite{Pa} 
\[
\frac 4{\pi ^2}\dint\limits_0^\ell \frac 1{(\ell -k)^2+1}\frac
1{k^2+1}dk=\frac 8{\pi ^2}\frac{\ell \arctan \left( \ell \right) +\ln \left(
\ell ^2+1\right) }{\left( \ell ^2+4\right) \ell } 
\]
for $\ell >0$. Why is this model incorrect? Clearly $\limfunc{E}(\theta
)=\pi /2$, not $\pi /3$, hence the three contact points are not equidistant
(on average). Consider also the triangle $T$ determined by the three lines: $%
\ell $ is the vertical side of $T$, but cannot be regarded as an
``arbitrary'' side. The assumption that $\theta \sim $ Uniform$[0,\pi ]$
requires change.

Another change (less a requirement than a preference) involves the
relationship between $T$ and the unit circle $C$. Clearly $C$ is an incircle
of $T$ if and only if there is no semicircle containing all three contact
points. Otherwise $C$ is an excircle of $T$. We wish to refine our model
(which presently incorporates both incircles and excircles) so that the
density of $\ell $ is based on incircles alone. Naturally $\ell >2$; the
infimum $2$ occurs in the limit as second and third lines both become
horizontal. The density for this refined model is given in the next section;
a related optimization problem appears at the end. As far as we know, these
results have not appeared in the random triangle literature before \cite{F1}.

\section{Unit Inradius}

Without loss of generality, let the first tangent be the vertical line
passing through the point $(-1,0)$. Let the second tangent pass through the
point $(\cos (\alpha ),\sin (\alpha ))$; let the third tangent pass through
the point $(\cos (\beta ),-\sin (\beta ))$. It is assumed that the bivariate
density for angles $\alpha $, $\beta $ is 
\[
\left\{ 
\begin{array}{lll}
2/\pi ^2 &  & \text{if }0<\alpha <\pi \text{, }0<\beta <\pi \text{ and }%
\alpha +\beta <\pi , \\ 
0 &  & \text{otherwise.}
\end{array}
\right. 
\]
It is best to think of $\alpha $ being measured in a counterclockwise
direction (as is customary) and $\beta $ being measured in a clockwise
direction. The condition $\alpha +\beta <\pi $ prevents contact points from
all crowding onto any semicircle (think of what happens when $\alpha +\beta
=\pi $). Dependency between $\alpha $ and $\beta $ makes our analysis more
complicated than earlier.

As a check, the univariate density for $\alpha $ is 
\[
\left\{ 
\begin{array}{lll}
2(\pi -\alpha )/\pi ^2 &  & \text{if }0<\alpha <\pi \text{,} \\ 
0 &  & \text{otherwise.}
\end{array}
\right. 
\]
Thus points on $C$ far away from $(-1,0)$ are favored (that is, small angles 
$\alpha $ are weighted more heavily than large $\alpha $) and $\limfunc{E}%
(\alpha )=\pi /3$.

We wish to determine the bivariate density of $h=\cot (\alpha /2)$, $k=\cot
(\beta /2)$. Via a Jacobian determinant argument, the density is \cite{Pa} 
\[
\left. \frac{\frac 2{\pi ^2}}{\frac 12\csc \left( \frac \alpha 2\right)
^2\frac 12\csc \left( \frac \beta 2\right) ^2}\right| \Sb \alpha =2\limfunc{%
arccot}(h),  \\ \beta =2\limfunc{arccot}(k)  \endSb =\frac 8{\pi ^2}\frac
1{h^2+1}\frac 1{k^2+1} 
\]
for $h>0$, $k>0$, $h\,k>1$. The latter inequality is true because 
\[
\pi =\alpha +\beta =2\limfunc{arccot}(h)+2\limfunc{arccot}(k) 
\]
if and only 
\[
k=\cot \left[ \frac 12\left( \pi -2\limfunc{arccot}(h)\right) \right] =\frac
1h. 
\]

Let $\ell =h+k$. The density of $\ell $ is the convolution 
\[
\frac 8{\pi ^2}\dint\limits_{a(\ell )}^{b(\ell )}\frac 1{(\ell -k)^2+1}\frac
1{k^2+1}dk 
\]
where limits of integration $a(\ell )$, $b(\ell )$ are found from $(\ell
-k)k>1$, hence $k^2-\ell \,k+1<0$; the zeroes of the quadratic are 
\[
\begin{array}{ccc}
a(\ell )=\dfrac{\ell -\sqrt{\ell ^2-4}}2>0, &  & b(\ell )=\dfrac{\ell +\sqrt{%
\ell ^2-4}}2<\ell
\end{array}
\]
and these are real because $\ell >2$. Integrating, we obtain the density of $%
\ell $ to be 
\[
\frac{16}{\pi ^2}\frac{f(\ell )+g(\ell )}{\left( \ell ^2+4\right) \ell } 
\]
for $\ell >2$, where 
\[
f(\ell )=\ell \arctan \left( \dfrac{\ell +\sqrt{\ell ^2-4}}2\right) -\ell
\arctan \left( \dfrac{\ell -\sqrt{\ell ^2-4}}2\right) , 
\]
\[
g(\ell )=\ln \left( \ell +\sqrt{\ell ^2-4}\right) -\ln \left( \ell -\sqrt{%
\ell ^2-4}\right) . 
\]
The mean of $\ell $ is infinite; its median $5.5482039188784452776442997...$
can be computed to high numerical precision as a consequence of our exact
density formula. See \cite{F2} for experimental confirmation of our work.

Having derived the density of an arbitrary side, let us briefly mention
other properties. The triangle $T$, under the condition that it
circumscribes the circle $C$, is acute with probability $1/4$ \cite{S1, S2}.
Since the inradius of $T$ is $1$, the area of $T$ is one-half the perimeter
of $T$ \cite{FO, OC} Unfortunately the perimeter density of $T$ (as well as
a trivariate density for sides) remains analytically intractible.

\section{Optimization Problem}

Of all triangles circumscribing the unit circle, an equilateral triangle
minimizes the perimeter $p$. The minimum value for $p$ is $6\sqrt{3}%
=10.39230...$.

Let $s$ denote a side of a triangle of unit inradius. If the other two sides
are nearly parallel and infinite, then $s$ approaches $2$ from above. The
infimum for $s$ is $2$.

Let $u$, $v$ denote two sides of a triangle of unit inradius. What is the
minimum value for $u+v$? It is surprising that this question is not better
known, especially since the answers for one side ($s$) and for the sum of
three sides ($p$) are clear.

By symmetry, the minimizing triangle is isosceles and $u=v$. Let $w$ denote
the remaining side of the triangle. By Heron's formula, 
\[
u+v+w=\frac 12\sqrt{(u+v+w)(-u+v+w)(u-v+w)(u+v-w)} 
\]
hence 
\[
4(u+v+w)=(-u+v+w)(u-v+w)(u+v-w) 
\]
hence 
\[
4(2v+w)=w^2(2v-w) 
\]
hence 
\[
v=\frac{\left( w^2+4\right) w}{2(w^2-4)}. 
\]
Differentiating with respect to $w$, we find that 
\[
w=\sqrt{8+4\sqrt{5}} 
\]
is a zero of the derivative. Substituting into the expression for $v$, we
deduce that the minimum value for $u+v$ is 
\[
2v=\sqrt{22+10\sqrt{5}}=6.66038.... 
\]

The angle $\theta $ at the apex of the minimizing triangle is also
interesting. By the Law of Cosines, 
\[
u^2+v^2-2u\,v\cos (\theta )=w^2
\]
hence 
\[
2v^2-w^2=2v^2\cos (\theta )
\]
hence 
\[
\cos (\theta )=1-\frac 12\frac{w^2}{v^2}=-2+\sqrt{5}=\frac 1{\varphi ^3}
\]
where $\varphi $ is the Golden mean \cite{F3}. Finally, $\theta
=1.33247...\approx 76.34^{\circ }$. This material constitutes a (very small)
first step toward characterizing the density for the sum of two arbitrary
sides of $T$.

\end{document}